\input amstex
\documentstyle{amsppt}
%\magnification=1200
\magnification=\magstephalf
\TagsOnRight
\NoBlackBoxes

\refstyle{A}
\widestnumber\key{BP2}

\topmatter
\NoBlackBoxes
\title Hyperbolic Reinhardt Domains in ${\Bbb C}^2$\\ 
with Noncompact Automorphism Group
\endtitle
\leftheadtext{\bf A.\ V.\ ISAEV AND S.\ G.\ KRANTZ}
\rightheadtext{\bf HYPERBOLIC REINHARDT DOMAINS}

%\footnote[]{{\bf Mathematics
%Subject Classification:} 32A07, 32H05, 32M05 \hfill}
%\footnote[]{{\bf Keywords and Phrases:} Automorphism
%groups, Reinhardt domains, hyperbolic. \hfill}
\author A. V. Isaev \ \ \ and \ \ \ S. G. Krantz
\endauthor

\address
Centre for Mathematics and Its Applications,
The Australian National University,
Canberra, ACT 0200,
AUSTRALIA 
\endaddress 
\email Alexander.Isaev\@anu.edu.au\endemail

\address
Department of Mathematics,
Washington University, St.~Louis, MO 63130,
USA 
\endaddress
\email sk\@math.wustl.edu
\endemail

\abstract {We give an explicit description of 
hyperbolic Reinhardt domains $D \subset {\Bbb C}^2$ 
such that: {\bf (i)}  $D$ has $C^k$-smooth boundary for
some $k \geq 1$,
{\bf (ii)}  $D$ intersects at least one of the coordinate
complex lines $\{z_1=0\}$, $\{z_2=0\}$, and {\bf (iii)}
$D$ has noncompact automorphism group.  We also give an
example that explains why such a setting is natural
for the case of hyperbolic domains and an example that
indicates that the situation in ${\Bbb C}^n$ for $n\ge
3$ is essentially more complicated than that in ${\Bbb
C}^2$.}  
\endabstract
\endtopmatter    
\document 

\heading 0. Introduction and Results\endheading

Let $D$ be a Kobayashi-hyperbolic
domain in ${\Bbb C}^n$, $n\ge 2$ (see \cite{Ko}
for terminology). Denote by $\text{Aut}(D)$ the
group of holomorphic automorphisms of $D$. The group
$\text{Aut}(D)$ with the topology of uniform convergence on
compact subsets of $D$ (the compact-open topology)
is in fact a Lie group (see
\cite{Ko}). 

The present paper is motivated by results characterizing a
domain by its automorphism group (see e.g. \cite{R}, \cite{W},
\cite{BP1}, \cite{BP2}). More precisely, we assume that
$\text{Aut}(D)$ is noncompact in the compact-open
topology. Most of the known results deal with the case
of bounded domains (see, however, \cite{B}). In the
present paper we consider possibly unbounded hyperbolic
domains.  Our thesis is that (unbounded) hyperbolic
domains have some of the geometric characteristics
of bounded domains.  In particular, they are tractable
for our studies.  But they also exhibit new automorphism
group action phenomena, and are therefore of special interest.
We present some of these new features in this work.

Here we assume that $D$ is a Reinhardt domain, i.e. a
domain which the standard action of the $n$-dimensional
torus ${\Bbb T}^n$ on ${\Bbb C}^n$,   
$$
z_j\mapsto
e^{{i\phi}_j}z_j,\qquad {\phi}_j\in
{\Bbb R},\quad j=1,\dots,n,\tag{1} 
$$
leaves invariant. In \cite{FIK} we gave a complete
classification of smoothly bounded Reinhardt domains
with noncompact automorphism group, and in \cite{IK} we
extended this result to Reinhardt domains with boundary
of any finite smoothness $C^k$, $k\ge 1$. One of the main
steps for obtaining these classifications was to
show that the noncompactness of $\text{Aut}(D)$ is
equivalent to that of $\text{Aut}_0(D)$, the connected
component of the identity in $\text{Aut}(D)$. We will now
explain this point in more detail, as it will provide some
motivation for the results of the present paper.

Following \cite{Sh}, we denote by
$\text{Aut}_{\text{alg}}(({\Bbb C}^*)^n)$ the group of
algebraic automorphisms of $({\Bbb C}^*)^n$, i.e. the
group of mappings of the form 
$$
z_i\mapsto\lambda_iz_1^{a_{i1}}\dots
z_n^{a_{in}},\quad i=1,\dots, n,\tag{2}
$$
where $\lambda_i\in{\Bbb C}^{*} \equiv {\Bbb C}\setminus\{0\}$,
$a_{ij}\in{\Bbb Z}$, and $\text{det}(a_{ij})=\pm 1$. For
a hyperbolic Reinhardt domain $D\subset{\Bbb C}^n$,
denote by $\text{Aut}_{\text{alg}}(D)$ the subgroup of
$\text{Aut}(D)$ that consists of algebraic
automorphisms  of $D$, i.e. automorphisms induced by
mappings from $\text{Aut}_{\text{alg}}(({\Bbb C}^*)^n)$.
It is shown in \cite{Kr} that
$\text{Aut}(D)=\text{Aut}_0(D)\cdot
\text{Aut}_{\text{alg}}(D)$, where the `` \, $\cdot$ \, '' denotes the
composition operation in $\text{Aut}(D)$. Therefore if
one can show that, for a hyperbolic Reinhardt domain
$D$, $\text{Aut}_{\text{alg}}(D)$ is finite up to the
action of ${\Bbb T}^n$ (see (1)), then the
noncompactness of $\text{Aut}(D)$ is equivalent to that
of $\text{Aut}_0(D)$ (see Proposition 1.1 in \cite{FIK}
for the case of bounded domains). Next, as is shown in
\cite{Kr}, $\text{Aut}_0(D)$ admits an explicit
description if $D$ is mapped into its normalized form by
a mapping of the form (2). This normalized form was the
main tool that we used in \cite{FIK}, \cite{IK}.

Unfortunately, as the following example shows, for the
case of hyperbolic Reinhardt domains the group
$\text{Aut}_{\text{alg}}(D)$ may be essentially
infinite, and therefore the scheme used in
\cite{FIK}, \cite{IK}, fails. 
\medskip

\noindent {\bf Example 1.} Consider the Reinhardt domain $D\subset
{\Bbb C}^2$ 
$$
D=\left\{\sin\left(\log\frac{|z_1|}{|z_2|}\right)<
\log |z_1
z_2|<\sin\left(\log\frac{|z_1|}{|z_2|}\right)+\frac{1}{2}\right\}.\tag{3}
$$
The boundary of $D$ is clearly $C^{\infty}$-smooth. The
group $\text{Aut}_{\text{alg}}(D)$ is not finite up to
the action of ${\Bbb T}^2$, since it contains
all the mappings 
$$
\aligned
&z_1\mapsto e^{\pi k}z_1,\\
&z_2\mapsto e^{-\pi k}z_2,
\endaligned
$$
for $k\in{\Bbb Z}$. This also shows, of course,
that $\text{Aut}(D)$ is noncompact.

To see that $D$ is hyperbolic, consider the mapping
$f\: D\rightarrow {\Bbb C}$, $f(z_1,z_2)=z_1 z_2$. It
is easy to see that $f$ maps $D$ onto the annulus
$A=\left\{e^{-1}<|z|<e^{\frac{3}{2}}\right\}$ which is a
hyperbolic domain in ${\Bbb C}$. The domains
$$
\aligned
&A_1=\left\{e^{-\frac{1}{4}}<|z|<e^{\frac{1}{2}}\right\},
\\ &A_2=\left\{e^{-1}<|z|<e^{-\frac{1}{8}}\right\}, \\
&A_3=\left\{e^{\frac{1}{4}}<|z|<e^{\frac{3}{2}}\right\}
\endaligned
$$
obviously cover $A$, and each of the preimages
$D_j=f^{-1}(A_j)$, $j=1,2,3$, is hyperbolic since
$D_j$ is contained in a union of bounded
pairwise nonintersecting domains. It then follows (see
\cite{PS}) that $D$ is hyperbolic.\qed
\medskip

It should be noted here that the domain (3)
does not intersect the coordinate complex lines
$\{z_1=0\}$, $\{z_2=0\}$ (note that, for any hyperbolic
Reinhardt domain in ${\Bbb C}^n$ not intersecting the
coordinate hyperplanes, one has
$\text{Aut}(D)=\text{Aut}_{\text{alg}}(D)$ \cite{Kr}).
As the following proposition shows, in complex dimension
$n=2$, the sort of pathology described in Example 1
above cannot occur if the domain intersects at least one
of the coordinate complex lines.

\proclaim{Proposition {\bf A}} Let $D\subset{\Bbb C}^2$ be
a hyperbolic Reinhardt domain with $C^1$-smooth
boundary, and let $D$ intersect at least one of the
coordinate complex lines $\{z_j=0\}$, $j=1,2$. Then
$\text{Aut}_{\text{alg}}(D)$ is finite up to the action
of ${\Bbb T}^2$.

In particular, for such a domain $D$, $\text{Aut}(D)$ is
noncompact if and only if $\text{Aut}_0(D)$ is
noncompact. 
\endproclaim

The above proposition allows us to use the description of
$\text{Aut}_0(D)$ from \cite{Kr} to obtain the following
classification result.

\proclaim{Theorem} Let $D\subset{\Bbb C}^2$ be
a hyperbolic Reinhardt domain with $C^k$-smooth
boundary, $k\ge 1$, and let $D$ intersect at least one of
the coordinate complex lines $\{z_j=0\}$, $j=1,2$. Assume
also that $\text{Aut}(D)$ is noncompact. Then $D$ is
biholomorphically equivalent to one of the following
domains: $$
\aligned
&\text{(i)}\,\left\{|z_1|^2+|z_2|^{\frac{1}{\alpha}}<1\right\},\\
&\qquad \text{where either $\alpha<0$, or
$\alpha=\frac{1}{2m}$ for some $m\in{\Bbb N}$, or 
$\alpha\ne\frac{1}{2m}$ for any}\\
& \qquad \text{$m\in{\Bbb N}$ and $0<\alpha<\frac{1}{2k}$;}\\
&\text{(ii)}\,\left\{|z_1|<1,(1-|z_1|^2)^{\alpha}<|z_2|<
R(1-|z_1|^2)^{\alpha}\right\},\\
& \qquad \text{where $1<R\le\infty$ and $\alpha<0$;}\\
&\text{(iii)}\,\left\{e^{\beta|z_1|^2}<|z_2|<Re^{\beta|z_1|^2}\right\},
\qquad \text{where $1<R\le\infty$, $\beta\in{\Bbb R}$, $\beta\ne 0$.}
\endaligned
$$

If $k<\infty$ and $\partial D$ is not
$C^{\infty}$-smooth, then $D$ is biholomorphically
equivalent to a domain of the form (i) with
$\alpha\ne\frac{1}{2m}$ for any
$m\in{\Bbb N}$ and $0<\alpha<\frac{1}{2k}$.

In case (i) the equivalence is given by dilations and
a permutation of the coordinates; in cases (ii) and
(iii) the equivalence is given by a mapping of the form
$$
\aligned
&z_{\sigma(1)}\mapsto \lambda z_{\sigma^{'}(1)}
z_{\sigma^{'}(2)}^a,\\ 
&z_{\sigma(2)}\mapsto \mu z_{\sigma^{'}(2)}^{\pm 1},
\endaligned
$$
where $\lambda,\mu \in {\Bbb C}^*$, $a\in{\Bbb Z}$ and
$\sigma$, $\sigma^{'}$ are permutations of $\{1,2\}$.
\endproclaim

In the following example we show that, in complex
dimension $n\ge 3$, an explicit classification result
analogous to the above theorem does not exist if we do
not impose extra conditions on the domain $D$. 
\medskip

{\bf Example 2.} Consider the domain $D\subset{\Bbb
C}^3$ given by
$$
\aligned
D&=\Biggl\{z : \phi(z) \equiv
|z_1|^2+(1-|z_1|^2)^2|z_2|^2\rho\Bigl(|z_2|^2(1-|z_1|^2),
|z_3|^2(1-|z_1|^2)\Bigr) \\ 
&\qquad \qquad \qquad \qquad
\qquad +(1-|z_1|^2)^2|z_3|^2 - 1 < 0 \Biggr\},
\endaligned\tag{4} $$
where $\rho(x_1,x_2)$ is a
$C^{\infty}$-smooth function on ${\Bbb R}^2$ such that
$\rho(x_1,x_2)>c>0$ everywhere, and the partial
derivatives of $\rho$ are  nonnegative for $x_1,x_2\ge 0$.

To show that $\partial D$ is smooth, we calculate
$$
\aligned
\frac{\partial\phi}{\partial
z_1}&=\overline{z_1}\Biggl(1-(1-|z_1|^2)\Bigl(2|z_2|^2\rho+
(1-|z_1|^2)|z_2|^4\frac{\partial\rho}{\partial
x_1}\\
& \qquad \qquad \qquad
+(1-|z_1|^2)|z_2|^2|z_3|^2\frac{\partial\rho}{\partial
x_2}+2|z_3|^2\Bigr)\Biggr),\\
\frac{\partial\phi}{\partial
z_2}&=(1-|z_1|^2)^2\overline{z_2}\left(\rho+(1-|z_1|^2)|z_2|^2
\frac{\partial \rho}{\partial x_1}\right),\\
\frac{\partial\phi}{\partial
z_3}&=(1-|z_1|^2)^2\overline{z_3}\left((1-|z_1|^2)|z_2|^2
\frac{\partial\rho}{\partial x_2}+1\right).
\endaligned\tag{5}
$$
It follows from (5) that not all the partial derivatives of
$\phi$ can vanish simultaneously at a point of
$\partial D$. Indeed, if $\frac{\partial\phi}{\partial
z_3}(p)=0$ at some point $p\in\partial D$ then, at $p$,
either $|z_1|=1$ or $z_3=0$. If $|z_1|=1$, then clearly
$\frac{\partial\phi}{\partial z_1}(p)\ne 0$. If $|z_1|\ne
1$, $z_3=0$, and, in addition,
$\frac{\partial\phi}{\partial z_2}(p)=0$, then $z_2=0$,
and therefore $|z_1|=1$, which is a contradiction.
Therefore, $\partial D$ is $C^{\infty}$-smooth.

To show that $D$ is hyperbolic, consider the holomorphic
mapping $f(z_1,z_2,z_3)=z_1$ from $D$ into ${\Bbb C}$.
Clearly, $f$ maps $D$ onto the unit disc
$\Delta=\{|z|<1\}$, which is a hyperbolic domain in
${\Bbb C}$. Further, the discs $\Delta_r=\{|z|<r\}$ for
$r<1$ form a cover of $\Delta$, and $f^{-1}(\Delta_r)$
is a bounded open subset of $D$ for any such $r$. Thus,
as in Example 1 above, we see that $D$ is
hyperbolic (see \cite{PS}).

Further, $\text{Aut}(D)$ is noncompact since it contains
the automorphisms
$$
(z_1, z_2, z_3) \longmapsto \left ( \frac{z_1-a}{1-\overline{a}z_1}\, , \,
   \frac{(1-\overline{a}z_1)z_2}{\sqrt{1-|a|^2}}\, , \,
   \frac{(1-\overline{a}z_1)z_3}{\sqrt{1-|a|^2}}  \right
)\tag{6} $$
for $|a|<1$.\qed 
\medskip

Examples similar to Example 2 can be constructed in
any complex dimension\linebreak $n\ge 3$. They indicate
that, most probably, there is no reasonable classification
of  smooth hyperbolic Reinhardt domains with noncompact
automorphism group for $n\ge 3$ even in the case when
the domains intersect the coordinate hyperplanes.
Indeed, in Example 2 we have substantial freedom in
choosing the function $\rho$. 

We note that the boundary of domain (4)
contains a complex hyperplane $z_1=\alpha$ for any
$|\alpha|=1$.  It may happen that, by imposing the
extra condition of the finiteness of type in the sense
of D'Angelo \cite{D'A} on the boundary of the domain, one
would eliminate the difficulty arising in Example 2 and
obtain an explicit classification. It also should be
observed that any point of the boundary of domain (4)
with $|z_1|=1$, $z_2=z_3=0$ is an orbit accumulation
point for $\text{Aut}(D)$ (see (6)); therefore, it is
plausible that one needs the finite type condition only
at such points (cf. the Greene/Krantz conjecture for
bounded domains \cite{GK}).

This work was initiated while the first 
author was an Alexander von Humboldt Fellow at the 
University of Wuppertal.  Research by the second author
was supported in part by NSF Grant DMS-9531967 and at MSRI
by NSF Grant DMS-9022140.

\heading Proof of Proposition A 
\endheading

We consider two cases.
\smallskip

{\bf Case 1.} Suppose first that $D$ intersects each of
the coordinate complex lines $\{z_j=0\}$, $j=1,2$. Then
any element $F\in\text{Aut}_{\text{alg}}(D)$ has the form
$$
\aligned
&z_1\mapsto \lambda z_{\sigma(1)}\\
&z_2\mapsto \mu z_{\sigma(2)}
\endaligned\tag{7}
$$
where $\lambda,\mu\in{\Bbb C}^*$ and
$\sigma$ is a permutation of the set $\{1,2\}$.

First let $\sigma=\text{id}$. We assume that
mapping (7) is not of the form (1);
hence either $|\lambda|\ne 1$, or $|\mu|\ne 1$. 
By passing to the
inverse mapping if necessary, we can also assume that
$|\lambda|<1$, or $|\mu|<1$.

Let $|\lambda|<1$. Take a point $p\in D$ of the form
$p=(c,0)$ and apply the $k^{\roman t \roman h}$ iteration
$F^k$ of $F$ to it: $F^k(p)=(\lambda^k c,0)$. Since
$|\lambda|<1$, it follows that $\lambda^k c\rightarrow 0$ as $k\rightarrow
\infty$, and therefore $(0,0)\in\overline{D}$. Since
$\partial D$ is $C^1$-smooth, we actually obtain that
$(0,0)\in D$. Therefore, for some $\epsilon>0$, the disc
$\Delta_{\epsilon}=\{|z_1|<\epsilon, z_2=0\}$ lies in
$D$. By applying the $k^{\roman t \roman h}$ iteration of $F^{-1}$ to
$\Delta_{\epsilon}$ and letting $k\rightarrow \infty$, we
obtain (since $|\lambda^{-k}|\rightarrow \infty$) that the
domain $D$ contains the entire complex line $\{z_2=0\}$
and therefore cannot be hyperbolic. The case of $|\mu|<1$
is treated similarly. Hence, $|\lambda|=|\mu|=1$, and $F$
is of the form (1).

Suppose now that $\sigma(1)=2$, $\sigma(2)=1$. We will
show that there exists no more than one automorphism of
the form (7) with this $\sigma$ (up to mappings of the
form (1)). Let $F_1$, $F_2$ be two such automorphisms,
with $F_j$ for $j=1,2$ given by
$$
\aligned
&z_1\mapsto \lambda_j z_2,\\
&z_2\mapsto \mu_j z_1,
\endaligned
$$
where $\lambda_j,\mu_j\in{\Bbb C}^*$.
Then, for the composition $F_1\circ F_2^{-1}$, we find that
$$
\aligned
&z_1\mapsto \frac{\lambda_1}{\lambda_2}z_1,\\
&z_2\mapsto \frac{\mu_1}{\mu_2}z_2.
\endaligned
$$
Hence, by the preceding argument,
$|\lambda_1|=|\lambda_2|$ and $|\mu_1|=|\mu_2|$; 
therefore $F_1$ differs from $F_2$ by a mapping of the
form (1).
\medskip

{\bf Case 2.} Let $D$ intersect only one of the
coordinate complex lines, say $\{z_1=0\}$. Then any
element of $\text{Aut}_{\text{alg}}(D)$ either has the
form
$$
\aligned
&z_1\mapsto \lambda z_1 z_2^a,\\
&z_2\mapsto \mu z_2,
\endaligned \tag{8}
$$
or the form
$$
\aligned
&z_1\mapsto \lambda z_1 z_2^a,\\
&z_2\mapsto \mu z_2^{-1},
\endaligned \tag{9}
$$
where $\lambda, \mu\in{\Bbb C}^*$, $a\in{\Bbb Z}$. We
will show that there is at most one element of
$\text{Aut}_{\text{alg}}(D)$ of each of the forms (8)
and (9).

Let $F_j$, $j=1,2$, be two automorphisms of the form
(8) given by
$$
\aligned
&z_1\mapsto \lambda_j z_1 z_2^{a_j},\\
&z_2\mapsto \mu_j z_2,
\endaligned
$$
where $\lambda_j,\mu_j\in{\Bbb C}^*$, $a_j\in{\Bbb Z}$.
Then for $F=F_1\circ F_2^{-1}$ we see that
$$
\aligned
&z_1\mapsto \frac{\lambda_1}
{\lambda_2}\mu_2^{a_2-a_1}z_1 z_2^{a_1-a_2},\\
&z_2\mapsto \frac{\mu_1}{\mu_2}z_2.
\endaligned
$$
Let $D_0=D\cap\{z_1=0\}$. Then, since $\partial D$ is
$C^1$-smooth, $\text{dist}(D_0,\{z_2=0\})>0$. Clearly,
$F$ preserves $D_0$. Suppose now that $|\mu_1|<|\mu_2|$.
Then, by considering the images of $D_0$ under iterations
of $F$, we see that $\text{dist}(D_0,\{z_2=0\})=0$ which
contradicts the smoothness of $\partial D$. Similarly, if 
$|\mu_1|>|\mu_2|$, then by applying iterations of $F^{-1}$
to $D_0$ we obtain the same contradiction. Therefore,
$|\mu_1|=|\mu_2|$.

By composing $F_2$ with a mapping of the form (1), we
can now assume that $\mu_1=\mu_2=\mu$ and therefore
$F$ is given by
$$
\aligned
&z_1\mapsto \frac{\lambda_1}
{\lambda_2}\mu^{a_2-a_1}z_1 z_2^{a_1-a_2},\\
&z_2\mapsto z_2.
\endaligned
$$
The $k^{\roman t \roman h}$ iteration of $F$ then has the form
$$
\aligned
&z_1\mapsto
\left(\frac{\lambda_1}{\lambda_2}\right)^k z_1
\left(\frac{z_2}{\mu}\right)^{k(a_1-a_2)} \\
&z_2\mapsto z_2.
\endaligned
$$
We now observe that there exist $\epsilon>0$ and a
disc $\tilde \Delta\subset{\Bbb C}$ such that the bidisc
$\{|z_1|<\epsilon, z_2\in\tilde\Delta\}$ lies in $D$.
Let $\Delta_{\epsilon,c}=\{|z_1|<\epsilon,z_2=c\}$ for
$c\in\tilde\Delta$. If for some $c\in
\tilde \Delta$ we have
$\left|\frac{\lambda_1}{\lambda_2}\right|\cdot 
\left|\frac{c}{\mu}\right|^{(a_1-a_2)}>1$ then, by
applying the iterations $F^k$ to $\Delta_{\epsilon,c}$
and letting $k\rightarrow \infty$, we see that the
complex line $\{z_2=c\}$ belongs entirely to $D$; this
conclusion contradicts the hyperbolicity of $D$. Similarly, if for
some $c\in \tilde\Delta$,
$\left|\frac{\lambda_1}{\lambda_2}\right|\cdot 
\left|\frac{c}{\mu}\right|^{(a_1-a_2)}<1$, then applying
iterations of $F^{-1}$ to $\Delta_{\epsilon,c}$ yields 
the same contradiction. Therefore,
$\left|\frac{\lambda_1}{\lambda_2}\right|\cdot 
\left|\frac{c}{\mu}\right|^{(a_1-a_2)}\equiv 1$ in
$\tilde\Delta$, and hence $a_1=a_2$,
$|\lambda_1|=|\lambda_2|$. Thus $F_1$ differs from
$F_2$ by a mapping of the form (1).

The case of mappings of the form (9) is treated
analogously. This completes the proof of
the proposition.\qed

\heading 2. Proof of Theorem\endheading

 We will use the following description of
$\text{Aut}_0(D)$ from \cite{Kr}. Any hyperbolic
Reinhardt domain in ${\Bbb C}^n$ can---by a
biholomorphic mapping of the form  (2)---be put into a
normalized form $G$ written as follows. There exist
integers $0\le s\le t\le p\le n$ and $n_i\ge 1$,
$i=1,\dots,p$, with $\sum_{i=1}^p n_i=n$, and
real numbers ${\alpha}_i^j$, $i=1,\dots,s$,
$j=t+1,\dots,p$, such that if we set
$z^i=\left(z_{n_1+\dots+n_{i-1}+1},\dots,z_{n_1+\dots+n_i}\right)$,
$i=1,\dots,p$, then $\tilde G:=G\bigcap\left\{z^i=0,\,
i=1,\dots,t\right\}$ is a hyperbolic Reinhardt domain in
${\Bbb C}^{n_{t+1}}\times\dots\times{\Bbb C}^{n_p}$, and
$G$ can be written in the form 
$$
\aligned
G=\Biggl\{\left|z^1\right|&<1,\dots,\left|z^s\right|<1,\\
&\Biggl(\frac{z^{t+1}}{\prod_{i=1}^s
\left(1-\left|z^i\right|^2\right)^{{\alpha}_i^{t+1}}\prod_{j=s+1}^t
\exp\left(-{\beta}_j^{t+1}\left|z^j\right|^2\right)}\ ,\ \dots \ , \\
&\frac{z^{p}}{\prod_{i=1}^s
\left(1-\left|z^i\right|^2\right)^{{\alpha}_i^p}\prod_{j=s+1}^t
\exp\left(-{\beta}_j^p\left|z^j\right|^2\right)}\Biggr)\in\tilde
G\Biggr\},
\endaligned\tag{10}
$$
for some real numbers ${\beta}_j^k$, $j=s+1,\dots,t$,
$k=t+1,\dots,p$. A normalized form can be chosen so that
$\text{Aut}_0(G)$ is given by the following formulas:
$$
\aligned
&z^i\mapsto\frac{A^iz^i+b^i}{c^iz^i+d^i},\quad i=1,\dots,s,\\
&z^j\mapsto B^jz^j+e^j,\quad j=s+1,\dots,t,\\
&z^k\mapsto
C^k\frac{\prod_{j=s+1}^t\exp\left(-\beta_j^k\left(2\overline{e^j}^TB^jz^j+
|e^j|^2\right)\right)z^k}{\prod_{i=1}^s(c^iz^i+d^i)^{2\alpha_i^k}},\quad
k=t+1,\dots,p,
\endaligned
$$
where
$$
\aligned
&\pmatrix
A^i&b^i\\
c^i&d^i
\endpmatrix\in SU(n_i,1),\quad i=1,\dots,s,\\
&B^j\in U(n_j),\quad e^j\in{\Bbb C}^{n_j},\quad j=s+1,\dots,t,\\
&C^k\in U(n_k),\quad k=t+1,\dots,p.
\endaligned
$$
The above classification implies that $\text{Aut}_0(G)$ is
noncompact only if $t>0$. 

Now let $n=2$. Clearly there are the following
possibilities  for a hyperbolic Reinhardt domain 
$\tilde G\subset{\Bbb C}$ (see (10)): 
\medskip 

\noindent  \text{\bf (i)} \ \ 
   $\displaystyle \tilde{G}=\{|z_2|<R\},\qquad 0<R<\infty;$
  
\noindent  \text{\bf (ii)} \   
   $\displaystyle \tilde{G}=\{r<|z_2|<R\},\qquad 0<r<R\le
\infty;$ 
     
\noindent  \text{\bf (iii)}    
   $\displaystyle
\tilde{G}=\{0<|z_2|<R\},\qquad
0<R<\infty.$  \medskip

This observation allows us to list all normalized forms
of hyperbolic Reinhardt domains in ${\Bbb C}^2$ as
follows
$$
\align
G&=\left\{|z_1|<1,|z_2|<R(1-|z_1|^2)^{\alpha}\right\},\\
& \qquad 0<R<\infty,\quad\alpha\in{\Bbb R},\tag{11}\\
G&=\left\{|z_1|<1,r(1-|z_1|^2)^{\alpha}<|z_2|<
R(1-|z_1|^2)^{\alpha}\right\},\\
&\qquad 0<r<R\le\infty,\quad\alpha\in{\Bbb R},\tag{12}\\
G&=\left\{|z_1|<1,0<|z_2|<
R(1-|z_1|^2)^{\alpha}\right\},\\
& \qquad 0<R<\infty,\quad\alpha\in{\Bbb R},\tag{13}\\
G&=\left\{re^{\beta|z_1|^2}<|z_2|<Re^{\beta|z_1|^2}\right\},\\
&\qquad 0<r<R\le\infty,\quad\beta\in{\Bbb R},\quad\beta\ne 0,\tag{14}\\ 
G&=\left\{0<|z_2|<Re^{\beta|z_1|^2}\right\},\\
&\qquad 0<R<\infty,\quad\beta\in{\Bbb R},\quad\beta\ne 0.\tag{15} 
\endalign 
$$
We are now going to select only those among the 
normalized forms (11)--(15) that can be the
images of domains with $C^k$-smooth boundary under
normalizing mappings of the form (2). 
We will treat each of cases (11)--(15)
separately.   
\medskip

{\bf Domain of type (11).} Observe first that, since the domain
$G$ contains the origin, the normalizing mapping
is linear (actually, it is given by dilations and a
permutation of the coordinates). Therefore, the
domain $G$ is a normalized form of a Reinhardt domain
with $C^k$-smooth boundary iff $\partial G$ is also
$C^k$-smooth. Hence $\alpha\ne 0$ (for otherwise $G$ is
a bidisc). If $\alpha>0$, then $G$ has a $C^k$-smooth
boundary iff either $\alpha=\frac{1}{2m}$,
for some $m\in{\Bbb N}$, or $\alpha\ne\frac{1}{2m}$
for any $m\in{\Bbb N}$ and $\alpha<\frac{1}{2k}$. If
$\alpha<0$, then $\partial G$ is $C^{\infty}$-smooth. It is
also clear that, for $k<\infty$, $\partial G$ has
$C^k$-smooth, but not $C^{\infty}$-smooth, boundary iff
$\alpha\ne\frac{1}{2m}$ for any $m\in{\Bbb N}$, and
$0<\alpha<\frac{1}{2k}$.     
\medskip

{\bf Domain of type (12).}  First of all, if $\alpha<0$,
then $\partial G$ is $C^{\infty}$-smooth. It is also clear
that, if $\alpha=0$, then $G$ cannot be a normalized form of
any Reinhardt domain with everywhere $C^k$-smooth
boundary for $k\ge 1$. 

Assume now that $\alpha>0$, and suppose first that
$R<\infty$. Then $\partial G$ is $C^{\infty}$-smooth
everywhere except at the points where $|z_1|=1$, $z_2=0$.
By applying the transformation 
$$
\aligned
&z_1\mapsto z_1,\\
&z_2\mapsto \frac{1}{z_2},
\endaligned\tag{16}
$$
we produce the following domain with $C^{\infty}$-smooth
boundary 
$$
\left\{|z_1|<1,\frac{1}{R}(1-|z_1|^2)^{-\alpha}
<|z_2|<\frac{1}{r}(1-|z_1|^2)^{-\alpha}\right\}.
$$

Let $\alpha>0$ and $R=\infty$. We claim that in this case
$G$ cannot be a normalized form of a Reinhardt domain
with everywhere $C^k$-smooth boundary for $k\ge 1$.
Indeed, this is easy to see if we notice that the
general form (up to permutation of the components) of
any mapping (2) that is biholomorphic on $G$ is as follows
$$
\aligned
&z_1\mapsto \lambda z_1 z_2^a,\\
&z_2\mapsto \mu z_2^{\pm 1},
\endaligned\tag{17}
$$
where $\lambda,\mu \in {\Bbb C}^*$, $a\in{\Bbb Z}$.

It is also easy to see that, for $k<\infty$, no domain
(12) can be a normalized form of a
Reinhardt domain with $C^k$-smooth, but not
$C^{\infty}$-smooth, boundary.
\medskip

{\bf Domain of type (13).} By transformation (16), $G$ is
mapped into a domain of the form (12) corresponding to
the case $R=\infty$, so it can be treated as above.
\medskip

{\bf Domain of type (14).} The boundary $\partial G$ of $G$ is
$C^{\infty}$-smooth. Also, if $k<\infty$, then $G$ cannot
be a normalized form of any Reinhardt domain with
$C^k$-smooth, but not $C^{\infty}$-smooth, boundary;
this assertion is proved by the same argument as we used for
domains of the form (12) above (see (17)).
\medskip

{\bf Domain of type (15).} By transformation (16), the domain $G$ is
mapped into a domain of the form (14) corresponding to
the case $R=\infty$, so it can be treated as above.
\medskip

The theorem is proved.\qed

\bigskip

\enddocument